\documentclass[11pt]{article}
\usepackage{amsthm}
\usepackage{latexsym, amssymb, amsmath}

\begin{document}
\numberwithin{equation}{section}

\newtheorem{THEOREM}{Theorem}
\newtheorem{PRO}{Proposition}
\newtheorem{XXXX}{\underline{Theorem}}
\newtheorem{CLAIM}{Claim}
\newtheorem{COR}{Corollary}
\newtheorem{LEMMA}{Lemma}
\newtheorem{REM}{Remark}
\newtheorem{EX}{Example}
\newenvironment{PROOF}{{\bf Proof}.}{{\ \vrule height7pt width4pt depth1pt} \par \vspace{2ex} }
\newcommand{\Bibitem}[1]{\bibitem{#1} \ifnum\thelabelflag=1 
  \marginpar{\vspace{0.6\baselineskip}\hspace{-1.08\textwidth}\fbox{\rm#1}}
  \fi}
\newcounter{labelflag} \setcounter{labelflag}{0}
\newcommand{\labelon}{\setcounter{labelflag}{1}}
\newcommand{\Label}[1]{\label{#1} \ifnum\thelabelflag=1 
  \ifmmode  \makebox[0in][l]{\qquad\fbox{\rm#1}}
  \else\marginpar{\vspace{0.7\baselineskip}\hspace{-1.15\textwidth}\fbox{\rm#1}}
  \fi \fi}

\newcommand{\LEFTLINE}{\ifhmode\newline\else\noindent\fi}
\newcommand{\RIGHTLINE}[1]{\LEFTLINE\rightline{#1}}
\newcommand{\CENTERLINE}[1]{\LEFTLINE\centerline{#1}}
\def\BOX #1 #2 {\framebox[#1in]{\parbox{#1in}{\vspace{#2in}}}}
\parskip=8pt plus 2pt
\def\AUTHOR#1{\author{#1} \maketitle}
\def\Title#1{\begin{center}  \Large\bf #1 \end{center}  \vskip 1ex }
\def\Author#1{\vspace*{-2ex}\begin{center} #1 \end{center}  
 \vskip 2ex \par}
\renewcommand{\theequation}{\arabic{section}.\arabic{equation}}
\def\bdk#1{\makebox[0pt][l]{#1}\hspace*{0.03ex}\makebox[0pt][l]{#1}\hspace*{0.03ex}\makebox[0pt][l]{#1}\hspace*{0.03ex}\makebox[0pt][l]{#1}\mbox{#1} }
\def\psbx#1 #2 {\mbox{\psfig{file=#1,height=#2}}}


\renewcommand{\theequation}{\arabic{equation}}
 
\Title{Nonnegative Trigonometric Polynomials \\ and Sturm's Theorem}

\begin{center}
MAN KAM KWONG\footnote{The research of this author is supported by the Hong Kong Government GRF Grant PolyU 5012/10P and the Hong Kong Polytechnic University Grants G-YK49 and G-U751}
\end{center}

\begin{center}
\emph{Department of Applied Mathematics\\ The Hong Kong Polytechnic University,\\ Hunghom, Hong Kong}\\
\tt{mankwong@polyu.edu.hk}\\[4ex]
July 2, 2015
\end{center}

\par\vspace*{\baselineskip}\par

\newcommand{\mb}{\mathbf}
\newcommand{\Cr}{\color{red}}

\begin{abstract}
\parskip=6pt
In an earlier article \cite{kw3}, we presented an algorithm
that can be used to rigorously check whether 
a specific cosine or sine polynomial is nonnegative in a given 
interval or not. The algorithm proves to be an indispensable
tool in establishing some recent results on nonnegative trigonometric polynomials. See, for example,
\cite{akw2}, \cite{kw1} and \cite{kw2}. It continues to play an essential
role in several ongoing projects.

The algorithm, however, cannot handle general trigonometric polynomials that involve
both cosine and sine terms. Some ad hoc methods to deal with
such polynomials have been suggested in \cite{kw3}, but none are, in general,
satisfactory. 

This note supplements \cite{kw3} by presenting an 
algorithm applicable to all general
trigonometric polynomials. It is based on the classical Sturm Theorem, 
just like the earlier algorithm. 

A couple of the references in \cite{kw3} are also updated.

\end{abstract}



\newpage

\section{A Summary of the Algorithm Proposed in \cite{kw3}}

\subsection*{\underline{$\vphantom{y}$Cosine Polynomials}}

\begin{enumerate}
\item Any cosine polynomial \,
$ C(x) = \sum_{k=1}^{n} \cos(kx)\, $
can be rewritten as an algebraic polynomial $ P(X) $ in the variable $ X=\cos(x) $.
The nonnegativity of $ C(x) $ in a given interval follows from that
of $ P(X) $ with $ X $ in a corresponding interval.
\item The nonnegativity of $ P(X) $ can be rigorously verified using the
classical Sturm Theorem (see \cite{vdw}).
\end{enumerate}

\noindent
With the help of the symbolic manipulation software MAPLE, these
steps can be automated with simple-to-use commands. The same is
true for the algorithm involving sine polynomials.

\subsection*{\underline{$\vphantom{y}$Sine Polynomials}}

\begin{enumerate}
\item Any sine polynomial \,
$ S(x) = \sum_{k=1}^{n} \sin(kx)\, $
can be rewritten as a product of $ \sin(x) $ and 
an algebraic polynomial $ P(X) $ in $ X $.
The nonnegativity of $ S(x) $ can be deduced from studying
the signs of $ P(X) $ in suitable subintervals of a corresponding interval.
\item The signs of $ P(X) $ can be determined rigorously using the
classical Sturm Theorem.
\end{enumerate}

\subsection*{\underline{$\vphantom{y}$General Trigonometric Polynomials}}
\quad  Some ad hoc methods are proposed, but none are satisfactory.

\section{New Algorithm for General Trigonometric Polynomials}

The new procedure is based on the elementary identities
\begin{equation}  \sin(x) = \frac{2T}{1+T^2} , \qquad  \cos(x) = \frac{1-T^2}{1+T^2} , \qquad  T = \tan \left( \frac{x}{2} \right)  \, .  \Label{scT}  \end{equation}

\begin{enumerate}
\item Any trigonometric polynomial 
can be rewritten as $ P_1(X)+\sin(x)\,P_2(X) $, where $ P_1(X) $ and $ P_2(X) $
are algebraic polynomials in $ X $.
\item Applying (\ref{scT}) results in a rational expression, the numerator
of which is an algebraic polynomial in $ T $ and the denominator is a power of 
$ (1+T^2) $.

The nonnegativity of the original trigonometric polynomial in a given interval follows from that
of the numerator of the rational expression in a corresponding interval.

\item The nonnegativity of the numerator can be rigorously verified using the
classical Sturm Theorem.
\end{enumerate}

\noindent
{\bf Remark.}
Although the new algorithm can also be used for a purely cosine polynomial or
a purely sine polynomial, it produces an algebraic polynomial with degree twice 
as that produced by the older algorithm. Hence, for purely cosine
or sine polynomials, the older algorithm is preferred.

\section{Examples and MAPLE Commands}
\begin{EX} \em
Let us illustrate the new procedure by proving that
\begin{equation}  A(x) = \frac{3}{5} + \sin(x) + \cos(x) + \frac{\sin(2x)}{2} + \frac{\cos(2x)}{2} \geq 0, \qquad  x \in [0,\pi ].  \Label{t1}  \end{equation}
$ A(x) $ can be rewritten as
$$  \frac{1}{10} +\cos \left( x \right) +\sin \left( x \right) + \cos^2 \left( x \right)  +\sin \left( x \right) \cos \left( x \right) .  $$
After applying (\ref{scT}), it becomes
\begin{equation}  {\frac {{T}^{4}-18\,{T}^{2}+40\,T+21}{ 10\left( 1+{T}^{2} \right) ^{2}}} \,.  \Label{t1y}  \end{equation}
As $ x $ varies from $ 0 $ to $ \pi  $, $ T $ varies from $ \tan(0) $ to $ \tan(\pi /2) $, or
from $ 0 $ to $ \infty  $.
Thus, (\ref{t1}) is true if we can show that the numerator of (\ref{t1y}) is nonnegative for
$ T\in[0,\infty ) $. The latter fact can be established by showing that the numerator,
which is a polynomial, has no root in $ [0,\infty ) $, using the classical Sturm Theorem.
\RIGHTLINE{\qed}
\end{EX}

\par\vspace*{\baselineskip}\par \hrule \par\vspace*{\baselineskip}\par

\subsection*{\underline{$\vphantom{y}$Using MAPLE}}
For more complicated trigonometric polynomials, the calculations required in the
two steps are often very
time-consuming and error-prone. They can be automated using MAPLE.

The following command defines a MAPLE function that takes a trigonometric polynomial as input and 
outputs the desired algebraic polynomial in $ T $.

\begin{quote}
\begin{verbatim}
pt := numer(subs(sin(x) = 2*T/(1+T^2),cos(x) = (1-T^2)/(1+T^2),
      expand(f)));
\end{verbatim}
\end{quote}

After issuing the above command, the next two affirm the assertion
of Example 1.

\begin{quote}
\begin{verbatim}
p := pt(3/5 + sin(x) + cos(x) + sin(2*x)/2 + cos(2*x)/2);
sturm(p, T, 0, infinity);
\end{verbatim}
\end{quote}

\noindent
The first command produces the numerator of the fraction in (\ref{t1y}) and assign
it to the variable {\tt p}. 
The second command uses the MAPLE built-in command {\tt sturm} to
determine the number of roots of $ p1 $ in
the specified interval $ T\in[0,\infty ) $.
In this example, the output from the second command is $ 0 $.

\par\vspace*{\baselineskip}\par
\par\vspace*{\baselineskip}\par \hrule \par\vspace*{\baselineskip}\par

\begin{EX} \em
Find the minimum of  
$ \displaystyle \sin(x) + \cos(x) + \frac{\sin(2x)}{2} + \frac{\cos(2x)}{2} , \quad  x\in[0,\pi ]. $

The problem is equivalent to the determination of the least value $ \alpha  $ such that
$$  B_\alpha (x) = \alpha  + \sin(x) + \cos(x) + \frac{\sin(2x)}{2} + \frac{\cos(2x)}{2} \geq  0 \quad  x\in[0,\pi ].  $$
Applying {\tt pt} to $ B_\alpha (x) $ produces the polynomial
$$  p_2(T) = (2\alpha -1)T^4 + (4\alpha -6)T^2 + 8T + (2\alpha +3) .  $$
Thus, we need to solve the equivalent problem of determining the least $ \alpha  $
such that $ p_2(T)\geq 0 $ in $ [0,\infty ) $. It is well-known that a necessary
condition for such an $ \alpha  $ is to be a root of the discriminant of $ p_2(T) $.
Using MAPLE, the discriminant is found to be
$$  \Delta  =  16384(2\alpha -1)(2\alpha -3)(8\alpha ^2+12\alpha -8).      $$
It is then easy to verify that one of the four roots of $ \Delta  $, namely,
$ \alpha =3(\sqrt{3}-1)/4 $, is the desired answer.

\qed
\end{EX}

\begin{EX} \em
Find the minimum of  
$ \displaystyle \sin(x) + \cos(x) + \sin(2x) + \cos(2x)  , \quad  x\in[0,\pi ]. $
The same arguments as in the previous example prove that the minimum
is a root of the equation (the lefthand side is the discriminant)
\begin{equation}  32768\,{a}^{4}-8192\,{a}^{3}-211968\,{a}^{2}+165888\,a+27648 = 0.  \end{equation}
However, the solutions no longer have simple representations. Numerical computation
gives the minimum as
$$  \alpha  \approx  1.040168473 .  $$

\qed
\end{EX}

\begin{EX} \em
In a recent project, we prove that
$$  \frac{\sin(x)}{3} + \frac{\sin(2x)}{2} + \sin(3x) + \frac{23}{125} \cdot 4 ,  $$
$$  \frac{\sin(x)}{4} + \frac{\sin(2x)}{3} + \frac{\sin(3x)}{2} + \sin(4x) +  \frac{23}{125} \cdot 5 ,  $$
(which are the first two cases of the more general polynomial)
$$  \sum_{k=1}^{n} \frac{\sin(kx)}{n-k+1} + \frac{23}{125} \cdot (n+1), \qquad  n = 3,4, \cdots ,  $$
are all nonnegative in $ [0,\infty ) $. Let us apply our new procedure to
the first two polynomials.

The corresponding algebraic polynomials are
$$  276\,{T}^{6}+1750\,{T}^{5}+828\,{T}^{4}-7000\,{T}^{3}+828\,{T}^{2}+ 3250\,T+276  $$
and
$$  138\,{T}^{8}-875\,{T}^{7}+552\,{T}^{6}+7375\,{T}^{5}+828\,{T}^{4}-9025 \,{T}^{3}+552\,{T}^{2}+1925\,T+138  $$
respectively. In both cases, the Sturm Theorem shows that the polynomials are nonnegative for 
$ T\in[0,\infty ) $.

\qed
\end{EX}

\end{document}